\documentclass[a4paper,12pt]{article}%
\usepackage{amsfonts}
\usepackage{amsmath}
\usepackage{amssymb}
\usepackage{graphicx}
\usepackage[onehalfspacing]{setspace}%
\newtheorem{theorem}{Theorem}[section]

\newtheorem{definition}{Definition}[section]

\newtheorem{lemma}{Lemma}[section]
\newtheorem{proposition}{Proposition}[section]
\newtheorem{remark}{Remark}[section]
\begin{document}

\title{The ARHD model}
\author{Andr\'{e} Mas$^{1}\thanks{Corresponding author : D\'{e}partement de
Math\'{e}matiques, CC 051, Universit\'{e} Montpellier 2, Place Eug\`{e}ne
Bataillon, 34095 Montpellier Cedex 5.\newline mas@math.univ-montp2.fr}$,
Besnik Pumo$^2$\\$^{(1)}${\small D\'{e}partment de Math\'{e}matiques, UMR\ 5149, Universit\'{e}
Montpellier 2,}\\$^{(2)}${\small Unit\'{e} de Statistiques, UMR A462 SAGAH, INH Angers}}
\date{}
\maketitle

\begin{abstract}
We introduce and study a new model for functional data. The ARHD is an
autoregressive model in which the first order derivative of the random curves
appears explicitely. Convergent estimates are obtained through an original
double penalization method. The prediction method is applied to a real set of
data already studied in the literature.

\end{abstract}


\section{Introduction\label{section.introduction}}

Usually time series may be viewed as the discretized observations obtained
from an underlying stochastic process $\left(  \xi(t),t\in\mathbb{R}\right)
$. Models and statistical inference on such processes aim at providing the
best possible predictor. Let us assume that the process $\xi$ is observed on
an interval $[0,T].$ We divide $[0,T]$ into $n$ subintervals $[i\delta
,(i+1)\delta]$, $i=0,\ldots,n-1$\ with $\delta=T/n$. This approach is clearly
justified in the case when $\xi$ is periodic with period $\delta$ but may be
generalized to processes that are stationary or not. In the following we
consider the functional-valued process $X=(X_{i},i\in\mathbb{Z})$ defined by:%
\[
X_{i+1}(t)=\xi(i\delta+t),\hspace{0.2in}0\leq t\leq\delta,\hspace{0.07in}%
i\in\mathbb{Z}.
\]
\ For a review on statistical analysis of functional data we refer to Ramsey
and Silverman (1997). In this paper we consider the prediction problem of the
process $\xi$\ on an entire time-interval $[T,T+\delta]$, or equivalently the
prediction of $X_{n+1}$\ knowing $X_{1},\ldots,X_{n}$. To deal with the
prediction problem Bosq (1991) introduced and studied an $H$\textit{-valued
autoregressive process of order one}, denoted ARH in the following, where
$\left(  H,\left\langle \cdot,\cdot\right\rangle \right)  $ is a suitable
Hilbert space of function with inner product $\left\langle \cdot
,\cdot\right\rangle $ such that $X_{i}\in H$ (typically $H=L^{2}\left[
0,\delta\right]  $ the space of square integrable functions on $\left[
0,\delta\right]  $). The model tends to generalize to functional data the
classical and celebrated $AR(1)$ model. Then the ARH process admits the
presentation
\[
X_{i}=\rho(X_{i-1})+\epsilon_{i}%
\]
\ where $(\epsilon_{i};i\in\mathbb{Z})$ \ is the $H$-valued innovation process
and $\rho$ is a bounded compact linear operator such that $\left\Vert
\rho\right\Vert _{\infty}<1.$ Note that $\left\Vert \rho\right\Vert _{\infty}$
stands for the classical operator norm for $\rho$ and is defined as%
\[
\left\Vert \rho\right\Vert _{\infty}=\sup_{h\in H,\left\Vert h\right\Vert
\leq1}\left\Vert \rho\left(  h\right)  \right\Vert .
\]
The ARH process is stationary and the best possible prediction of $X_{n+1}$ is
$\hat{X}_{n+1}=\rho(X_{n})$ whenever $\mathbb{E}\left(  \epsilon_{n}%
|X_{n-1}\right)  =0$. If $\widehat{\rho}_{n}$ is a consistent estimate for
$\rho,$ the prediction is made through $\widehat{\rho}_{n}(X_{n})$. Bosq
(2000) proposes such a predictor and proves the consistency under mild conditions.

Considering some regularity conditions on the sample paths, one may obtain
similar results for autoregressive spaces with values in other functional
spaces (see Pumo (1992) for results on $C[0,1]$, or Mourid (1995) for results
on general Banach spaces). Alternative approaches to solve the prediction
problem based on ARH modelization are proposed by Besse et al. (1996, 2000))
by means of \textit{spline smoothing} and more recently Antoniadis and
Sapatinas (2003) who implemented \textit{Wavelet techniques}. Furthermore
numerical studies show that the predictors obtained by this alternative
methods are better then those obtained by \textit{linear interpolation ARH
predictor} (Pumo (1998)).

Generalizing to the functional data the classical multivariate models is a new
and fruitful trend in modern statistics : linear or non linear regression,
high-dimensional or functional ANOVA, MA($\infty$) processes all have their
"functional" counterparts. But usually and up to the authors' knowledge, these
classical statistical models for functional data never involve the derivatives
of the random curves rebuilt from discretized data. However the opportunity to
compute explicit first (or higher) order derivatives is one of the main
feature differencing truly functional from multivariate data. We refer to
Silverman (1996) in the framework of principal component analysis, then
Ferraty and Vieu (2003) for regression models. These authors underline the
specific amount of information contained in the derivatives of curves rebuilt
from functional data as well as their practical interests.

The aim of this paper is to use the functional properties of smooth sample
paths in order to improve the predictor. More precisely we suppose that the
sample paths belongs to the Sobolev space $W^{2,1}$\textit{\ }(defined in the
next section). We study an autoregressive model, initially introduced by
Marion and Pumo (2004), whose definition explicitly involves the first order
derivative of the data. This "ARHD" model is detailed below. Wong process,
that is a stationary Gaussian process with continuously differentiable paths
may be represented this way. We propose to estimate the two unknown parameters
of the model by an original method. This technique is inspired by the ridge
regression method and involves two overlapping penalization through two
parameters depending on each other.

The paper is organized as follows. In the following paragraph we introduce the
ARHD(1) model, simply denoted ARHD in the sequel, and show that it is strictly
stationary. In section \ref{section.estimation.procedure} we give conditions
for the unknown operators $\varphi$ and $\Psi$ to be identifiable and provide
estimates as well as asymptotic results. Section \ref{section.applications} is
devoted to technical details about the numerical calculation of ARHD
predictors and to the comparison of the ARHD predictors with various
functional methods in two cases: the first one is a simulated example (the
Wong process which allows an ARHD presentation) ; the second is a real data
study concerning the El Ni\^{n}o-Southern Oscillation (ENSO) time series.
Proofs of asymptotic results are postponed to section \ref{section.proofs}.

\section{The model\label{sub.section.ARHD.model}}

It was introduced by Pumo and Marion (2004). Let $X_{i}$ be a sample of random
curves. We introduced the model above :%

\begin{equation}
X_{i+1}=\phi\left(  X_{i}\right)  +\Psi\left(  X_{i}^{\prime}\right)
+\varepsilon_{i+1} \label{ARHD}%
\end{equation}
where $\phi$ and $\Psi$ are linear operators.

Now we suppose that for all $i,$ $X_{i}$ takes its values in the Sobolev space
$W^{2,1}\left[  0,1\right]  .$%
\[
W^{2,1}=\left\{  u\in L^{2}\left[  0,1\right]  ,u^{\prime}\in L^{2}\left[
0,1\right]  \right\}  .
\]
The space $W^{2,1}$ is a separable Hilbert space endowed with scalar product :%
\[
\left\langle u,v\right\rangle _{W}=\int_{0}^{1}u\left(  t\right)  v\left(
t\right)  dt+\int_{0}^{1}u^{\prime}\left(  t\right)  v^{\prime}\left(
t\right)  dt.
\]
We refer to Ziemer (1989) or to Adams and Fournier (2003) for monographs
dedicated to Sobolev spaces. In the sequel $W^{2,1}$ will be denoted $W$ and
$W^{2,0}=L^{2}$ will be denoted $L$ for the sake of simplicity. Obviously if
we set $Du=u^{\prime}$ then $D$ maps $W$ onto $L$ ($D$ is the ordinary
differential operator)$.$ Furthermore Sobolev's imbedding theorem ensures that
(see Adams and Fournier (2003) Theorem 4.12 p.85)%
\[
\left\Vert Du\right\Vert _{L}\leq C\left\Vert u\right\Vert _{W}%
\]
(where $C$ is some constant which does not depend on $u$) i.e. $D$ is a
bounded operator from $W$ to $L.$

From now on we assume that $\phi$\ is\ a\ compact\ operator\ from\ $W$%
\ to\ $W$\ and\ $\Psi$\ is\ a\ compact\ operator\ from\ $L$\ to\ $W.$ For a
review on compact operators we refer to Dunford-Schwartz (1988) or Gohberg,
Goldberg, Kaashoek (1991).

\section{ARW representation of the ARHD process.}

From the above paragraph we know that $\phi+\Psi D$ is a well defined operator
on $W$ and onto $W,$ that $\Psi D$ is a compact operator as the product of a
bounded and a compact operator and consequently that $\phi+\Psi D$ is itself
compact as the sum of two compact operators.

We can rewrite (\ref{ARHD}) :%

\begin{equation}
X_{i+1}=A\left(  X_{i}\right)  +\varepsilon_{i+1} \label{ARH1}%
\end{equation}
where $A=\phi+\Psi D.$

Finally the ARHD process may be rewritten as a special ARH(1) process with
values in $W$.\bigskip

The trouble with (\ref{ARH1}) is the following : the parameters $\phi$ and
$\Psi$ are hidden behind $A$ and we are not willing to infer on the latter.
Obviously we are going to face two issues :

\begin{itemize}
\item Studying the identifiability of $\phi$ and $\Psi$ in the model above.

\item Providing a consistent estimation procedure for $\phi$ and $\Psi$ before forecasting.
\end{itemize}

From now on we suppose that%
\begin{align*}
\mathbf{H1}  &  :\left\Vert A\right\Vert _{\infty}<1,\\
\mathbf{H2}  &  :\left\Vert X\right\Vert _{W}<+\infty\quad a.s.
\end{align*}
The first assumption is crucial for the stationarity of the process. The
second is quite restrictive but could be alleviated to mild moment assumptions
but it will make the proofs of the main result more easily readable. This
assumption appears for instance in Cardot, Ferraty, Sarda (1999) for the same
reasons. We assess the first property of the process, which will be useful in
the sequel.

\begin{proposition}
When assumptions $\mathbf{H1}$ and $\mathbf{H2}$ hold, $\left(  X_{i}\right)
_{i\in\mathbb{Z}}$ and $\left(  X_{i}^{\prime}\right)  _{i\in\mathbb{Z}}$ are
strictly stationary sequences on $W$ and $L$ respectively.
\end{proposition}

The stationarity of $\left(  X_{i}\right)  _{i\in\mathbb{Z}}$ is a simple
consequence of the representation of equation (\ref{ARH1}) and of previous
results obtained for instance by Bosq (2000), Chapter 3. The continuity of $D$
on $W$ implies the stationarity of the sequence $\left(  X_{i}^{\prime
}\right)  _{i\in\mathbb{Z}}.$

\section{Estimation procedure\label{section.estimation.procedure}}

\subsection{The moment method}

From a practical point of view the Sobolev setting is not really a constraint.
It is well known that either splines or wavelets will provide standard
reconstruction method (from the discretized data) yielding functions in
$W^{2,1}.$

The model is purely functional : we cannot invoke any likelihood-based
technique since the "density of a random curve" makes non sense (Lebesgue's
measure does not exist on infinite dimensional spaces). We propose to start
from a classical moment method and to adapt it to our setting.

By $\mathcal{C}_{H}$ (resp. $\mathcal{C}_{HH^{\prime}}$) we denote the space
of compact operators on the Hilbert space $H$ (resp. mapping the Hilbert space
$H$ onto $H^{\prime}$). Some finite rank operators are defined by means of the
tensor product : if $u$ and $v$ belong to $H$ and $H^{\prime}$ respectively
$u\otimes_{H}v$ is the operator defined on $H$ by : for all $h\in H$%
\[
\left(  u\otimes_{H}v\right)  \left(  h\right)  =\left\langle u,h\right\rangle
_{H}v.
\]

We start from a sample $\left(  X_{i},X_{i}^{\prime}\right)  _{1\leq i\leq n}$
and we denote%
\begin{align*}
\Gamma &  =E\left(  X_{0}\otimes_{W}X_{0}\right)  ,\;\Gamma^{\prime}=E\left(
X_{0}\otimes_{W}X_{0}^{\prime}\right)  ,\\
\Gamma^{\prime\ast}  &  =E\left(  X_{0}^{\prime}\otimes_{L}X_{0}\right)
,\;\Gamma^{\prime\prime}=E\left(  X_{0}^{\prime}\otimes_{L}X_{0}^{\prime
}\right)  ,\\
\Delta &  =E\left(  X_{0}\otimes_{W}X_{1}\right)  ,\;\Delta^{\prime}=E\left(
X_{0}^{\prime}\otimes_{L}X_{1}\right)  .
\end{align*}
Under assumption $\mathbf{H2}$ all these operators belong either to
$\mathcal{C}_{W},$ $\mathcal{C}_{WL},$ $\mathcal{C}_{LW}$ or $\mathcal{C}_{L}%
$. In fact assumption $\mathbf{H2}$ could be replaced by $E\left\Vert
X\right\Vert _{W}^{2}<+\infty.$

By $\Gamma_{n},\Gamma_{n}^{\prime},...,\Delta_{n}^{\prime}$ we denote the
empirical counterpart of these operators based on the sample $\left(
X_{i},X_{i}^{\prime}\right)  _{1\leq i\leq n}$. For example :%
\begin{align}
\Gamma_{n}  &  =\dfrac{1}{n}\sum_{k=1}^{n}X_{k}\otimes_{W}X_{k},\label{esti}\\
\Delta_{n}^{\prime}  &  =\dfrac{1}{n-1}\sum_{k=1}^{n-1}X_{k}^{\prime}%
\otimes_{L}X_{k+1}.\nonumber
\end{align}

\begin{remark}
The notation $\Gamma^{\prime\ast}$ is not ambiguous : $\Gamma^{\prime\ast}$ is
truly the adjoint operator of $\Gamma^{\prime}$.
\end{remark}

\begin{remark}
Conversely, if the random function $X^{\prime}$ is truly the derivative of
$X,$ this is no more the case as far as linear operators are concerned :
$\Gamma^{\prime}$ is not the derivative of $\Gamma.$ The $^{\prime}$ is just a
notation in this setting ; this would make no sense anyway. However, it should
be remarked that for all $u$ in $L$%
\[
\left(  \Gamma^{\prime\ast}\left(  u\right)  \right)  ^{\prime}=\left(
E\left(  \left\langle u,X^{\prime}\right\rangle _{L}X\right)  \right)
^{\prime}=\Gamma^{\prime\prime}\left(  u\right)
\]
and for all $v$ in $W$%
\[
\left(  \Gamma\left(  v\right)  \right)  ^{\prime}=\left(  E\left(
\left\langle v,X\right\rangle _{W}X\right)  \right)  ^{\prime}=\Gamma^{\prime
}\left(  v\right)  .
\]

\end{remark}

Quite naturally, from (\ref{ARHD}) -multiply with $\left\langle X_{i}%
,\cdot\right\rangle $ and $\left\langle X_{i}^{\prime},\cdot\right\rangle $
successively then take expectation- we easily deduce both moment equations :%

\begin{equation}
\mathcal{S}=\left\{
\begin{tabular}
[c]{l}%
$\Delta=\phi\Gamma+\Psi\Gamma^{\prime}$\\
$\Delta^{\prime}=\phi\Gamma^{\prime\ast}+\Psi\Gamma^{\prime\prime}.$%
\end{tabular}
\ \ \right.  \label{Mom}%
\end{equation}
Resolving this system is apparently easy but we should be aware of two facts :

\begin{itemize}
\item Operators (here, $\Delta,\phi,\Gamma..$.) do not commute !

\item The inverse operators of $\Gamma,$ and $\Gamma^{\prime\prime}$ do not
necessarily exist and when they do, they are unbounded, i.e. not continuous
(remind that $\Gamma,$ and $\Gamma^{\prime\prime}$ are compact operators and
that compact operators have no bounded inverses).
\end{itemize}

At this point, before trying to solve (\ref{Mom}) we need to study
identifiability of the unknown infinite dimensional parameter $\left(
\phi,\Psi\right)  \in\mathcal{C}_{W}\times\mathcal{C}_{LW}$ in our statistical problem.

\subsection{Identifiability}

We set $\mathcal{E}=\mathcal{C}_{W}\times\mathcal{C}_{LW}.$ If both equations
in (\ref{Mom}) are the starting point we should make sure that solutions to
these equations are well and uniquely defined. Suppose for instance that
\textrm{Ker}$\Gamma\neq\left\{  0\right\}  $ and take $h$ in it. Now set
$\widetilde{\phi}=\phi+h\otimes_{W}h.$ Then
\[
\widetilde{\phi}\Gamma=\phi\Gamma+\left(  h\otimes_{W}h\right)  \Gamma
\]
but $\left(  h\otimes_{W}h\right)  \Gamma=0$. So $\widetilde{\phi}\Gamma
=\phi\Gamma$ and $\phi$ is not unique (there are even infinitely many
solutions in the space $\phi+$\textrm{Ker}$\Gamma$). The next assumption is%
\[
\mathbf{H3}:\mathrm{Ker}\Gamma=\mathrm{Ker}\Gamma^{\prime\prime}=\left\{
0\right\}  \mathrm{.}%
\]
In other words we suppose that both operators above are one to one.

Now turning back to (\ref{Mom}) we rewrite the system. Equivalently :%
\begin{align}
\mathcal{S}  &  =\left\{
\begin{tabular}
[c]{l}%
$\Delta=\phi\Gamma+\Psi D\Gamma$\\
$\Delta^{\prime}=\phi\Gamma D^{\ast}+\Psi D\Gamma D^{\ast}$%
\end{tabular}
\ \ \ \right. \label{Mom2}\\
&  \iff\left\{
\begin{tabular}
[c]{l}%
$\left(  \Delta,\Delta^{\prime}\right)  =\left(  \phi,\Psi\right)  \Lambda$\\
$\Lambda=\left(
\begin{array}
[c]{cc}%
\Gamma & \Gamma D^{\ast}\\
D\Gamma & D\Gamma D^{\ast}%
\end{array}
\right)  .$%
\end{tabular}
\ \ \ \right.  \label{mom3}%
\end{align}
We are now ready to solve the identification problem.

\begin{proposition}
\label{Ident}The couple $\left(  \phi,\Psi\right)  \in\mathcal{E}$ is
identifiable for the moment method proposed in (\ref{Mom}) if and only if
$\left(  \phi,\Psi\right)  \notin\mathcal{N}$ where $\mathcal{N}$ is the
vector subspace of $\mathcal{E}$ defined by%
\[
\mathcal{N=}\left\{  \left(  U,V\right)  \in\mathcal{E}:U+VD=0\right\}  .
\]
Note that $\mathcal{N}$ is a closed set in $\mathcal{E}$.
\end{proposition}

In other words if, for all $u$ in $W,$ $\phi u+\Psi u^{\prime}=0,$ the
parameter $\left(  \phi,\Psi\right)  $ cannot be identified.

The Proposition is proved at the beginning of the last section of the paper.

\section{Definition of the estimates and convergence}

Since the unknown parameters $\phi$ and $\Psi$ are operators estimating them
means dealing with random operators based on the double sample $\left(
X_{i},X_{i}^{\prime}\right)  _{1\leq i\leq n}$. We refer for instance to
(\ref{esti}) above for examples of these available operators.

The estimates stem from (\ref{Mom}) which is a highly non invertible system.
We are classically going to add a small perturbation to regularize it and make
it invertible. We solve :%

\[
\mathcal{S}^{\prime}=\left\{
\begin{tabular}
[c]{l}%
$\Delta=\phi\left(  \Gamma+\alpha I_{W}\right)  +\Psi\Gamma^{\prime}$\\
$\Delta^{\prime}=\phi\Gamma^{\prime\ast}+\Psi\left(  \Gamma^{\prime\prime
}+\alpha I_{L}\right)  $%
\end{tabular}
\ \ \ \ \ \ \right.
\]
where $\alpha$ is a positive real number and $I_{W}$ denotes the identity
operator on $W.$ Now the operators $\left(  \Gamma+\alpha I_{W}\right)  $ and
$\left(  \Gamma^{\prime\prime}+\alpha I_{L}\right)  $ are no more compact but
have bounded inverses. Basic algebra gives :%
\begin{equation}
\mathcal{S}^{\prime}=\left\{
\begin{tabular}
[c]{l}%
$\phi\left[  \left(  \Gamma+\alpha I_{W}\right)  -\Gamma^{\prime\ast}\left(
\Gamma^{\prime\prime}+\alpha I_{L}\right)  ^{-1}\Gamma^{\prime}\right]
=\Delta-\Delta^{\prime}\left(  \Gamma^{\prime\prime}+\alpha I_{L}\right)
^{-1}\Gamma^{\prime}$\\
$\Psi\left[  \left(  \Gamma^{\prime\prime}+\alpha I_{L}\right)  -\Gamma
^{\prime}\left(  \Gamma+\alpha I_{W}\right)  ^{-1}\Gamma^{\prime\ast}\right]
=\Delta^{\prime}-\Delta\left(  \Gamma+\alpha I_{W}\right)  ^{-1}\Gamma
^{\prime\ast}.$%
\end{tabular}
\ \ \ \ \ \ \ \ \ \right.  \label{sol}%
\end{equation}
Which is then once more approximated by :%
\[
\mathcal{S}^{\prime\prime}=\left\{
\begin{tabular}
[c]{l}%
$\phi\left[  \Gamma-\Gamma^{\prime\ast}\left(  \Gamma^{\prime\prime}+\alpha
I_{L}\right)  ^{-1}\Gamma^{\prime}\right]  =\Delta-\Delta^{\prime}\left(
\Gamma^{\prime\prime}+\alpha I_{L}\right)  ^{-1}\Gamma^{\prime}$\\
$\Psi\left[  \Gamma^{\prime\prime}-\Gamma^{\prime}\left(  \Gamma+\alpha
I_{W}\right)  ^{-1}\Gamma^{\prime\ast}\right]  =\Delta^{\prime}-\Delta\left(
\Gamma+\alpha I_{W}\right)  ^{-1}\Gamma^{\prime\ast}.$%
\end{tabular}
\ \ \ \ \ \ \ \ \ \right.
\]
We just dropped $\alpha I_{W}$ on the first line and $\alpha I_{L}$ on the
second to get $\mathcal{S}^{\prime\prime}$. Take the first line in the above
display. The operator
\begin{equation}
S_{\phi}=\Gamma-\Gamma^{\prime\ast}\left(  \Gamma^{\prime\prime}+\alpha
I_{L}\right)  ^{-1}\Gamma^{\prime} \label{Sphi}%
\end{equation}
is a selfadjoint compact operator. (Indeed $\Gamma^{\prime\ast}\left(
\Gamma^{\prime\prime}+\alpha I_{L}\right)  ^{-1}\Gamma^{\prime}$ is a compact
operator because $\Gamma^{\prime}$ and $\Gamma^{\prime\ast}$ are). We may
deduce from this fact that $S_{\phi}$ has real eigenvalues (not necessarily
positive) and furthermore that once again $S_{\phi}$ has no bounded inverse.
The same remarks hold for
\begin{equation}
S_{\Psi}=\Gamma^{\prime\prime}-\Gamma^{\prime}\left(  \Gamma+\alpha
I_{W}\right)  ^{-1}\Gamma^{\prime\ast}. \label{Spsi}%
\end{equation}
However we can provide an approximate solution to $\mathcal{S}^{\prime\prime}$
by regularizing $S_{\phi}$ and $S_{\Psi}$ once more by a penalization method.
Finally the pseudo solutions we propose to solve $\mathcal{S}^{\prime\prime}$
hence $\mathcal{S}$ are based on a second strictly positive parameter $\beta$
and are denoted $\widetilde{\phi}$ and $\widetilde{\Psi}$ :%
\begin{equation}
\left\{
\begin{tabular}
[c]{l}%
$\widetilde{\phi}=\left[  \Delta-\Delta^{\prime}\left(  \Gamma^{\prime\prime
}+\alpha I_{L}\right)  ^{-1}\Gamma^{\prime}\right]  \left(  S_{\phi}+\beta
I\right)  ^{-1}$\\
$\widetilde{\Psi}=\left[  \Delta^{\prime}-\Delta\left(  \Gamma+\alpha
I_{W}\right)  ^{-1}\Gamma^{\prime\ast}\right]  \left(  S_{\phi}+\beta
I\right)  ^{-1}.$%
\end{tabular}
\ \ \ \ \ \ \ \ \ \right.  \label{appro}%
\end{equation}

This new system defines relations from which we propose to deduce estimates.
From now on -in order to alleviate the notations- by $S^{\dagger}$ we denote
the operator defined by $\left(  S+\alpha_{n}I_{L}\right)  ^{-1}$ where
$\alpha_{n}$ is a non increasing sequence of positive numbers decaying to
zero. We set :%
\begin{align}
S_{n,\phi}  &  =\Gamma_{n}-\Gamma_{n}^{\prime\ast}\left(  \Gamma_{n}%
^{\prime\prime\dagger}\right)  \Gamma_{n}^{\prime},\label{Snphi}\\
S_{n,\Psi}  &  =\Gamma_{n}^{\prime\prime}-\Gamma_{n}^{\prime}\left(
\Gamma_{n}^{\dagger}\right)  \Gamma_{n}^{\prime\ast},\label{Snpsi}\\
T_{n,\phi}  &  =\Delta_{n}-\Delta_{n}^{\prime}\left(  \Gamma_{n}^{\prime
\prime\dagger}\right)  \Gamma_{n}^{\prime},\label{Tnphi}\\
T_{n,\Psi}  &  =\Delta_{n}^{\prime}-\Delta_{n}\left(  \Gamma_{n}^{\dagger
}\right)  \Gamma_{n}^{\prime\ast}. \label{Tnpsi}%
\end{align}
Taking $\beta_{n}\downarrow0$ we obtain the following

\begin{definition}
The estimate of the couple $\left(  \phi,\Psi\right)  $ is $\left(  \phi
_{n},\Psi_{n}\right)  $ based on (\ref{appro}) and defined by :%
\begin{equation}
\left\{
\begin{tabular}
[c]{l}%
$\phi_{n}=T_{n,\phi}\left(  S_{n,\phi}+\beta_{n}I\right)  ^{-1}$\\
$\Psi_{n}=T_{n,\Psi}\left(  S_{n,\phi}+\beta_{n}I\right)  ^{-1}.$%
\end{tabular}
\ \ \ \right.  \label{estim}%
\end{equation}

\end{definition}

The next Theorem is the main theoretical result of this article. It provides
the convergence of our estimates when the sample size goes to infinity.

\begin{theorem}
\label{Conv}When $\mathbf{H1-3}$ hold and if $\alpha_{n}\rightarrow0$,
$\beta_{n}\rightarrow0$ with $\sqrt{n}\alpha_{n}^{2}\beta_{n}^{2}%
\rightarrow+\infty$ and $\sqrt{\alpha_{n}}/\beta_{n}\rightarrow0,$%
\begin{align*}
&  \phi_{n}\overset{\mathbb{P}}{\rightarrow}\phi,\\
&  \Psi_{n}\overset{\mathbb{P}}{\rightarrow}\Psi.
\end{align*}

\end{theorem}

The convergence is understood in the $\left\Vert \cdot\right\Vert _{\infty}$
norm for bounded operators.

Note that Theorem \ref{Conv} holds whenever $\alpha_{n}=n^{-a}$ and $\beta
_{n}=n^{-b}$ with $b<a/2$ and $2b+2a<1/2.$

\begin{remark}
Originally the linear model (\ref{ARHD}) is subject to serious
multicollinearity troubles since $X_{n}^{\prime}=DX_{n}.$ Even if the curve
$X_{n}^{\prime}$ usually looks quite different from $X_{n},$ there is a total
stochastic dependence between them. The method used in this article to tackle
this problem (as well as the intrinsic "inverse problem" aspects related to
the inversion of the covariance operators $\Gamma$ and $\Gamma^{\prime\prime}%
$) is new up to the authors' knowledge. As it can be seen through above at
display (\ref{appro}) or in the proofs below, it relies on a double
penalization technique first by the index $\alpha_{n}$ then by $\beta_{n}$
linking both indexes in order to asymptotically suppress the bias terms.
\end{remark}

\section{A numerical study and application: ENSO \label{section.applications}}

\bigskip In this section we illustrate the ARHD method of prediction proposed
in this paper by some numerical studies for two examples. We give first some
technical results to carry out numerical calculations. The first application
is connected to Wong's process (see Wong (1966)) which admits an ARHD
presentation. We compare the ARHD predictor with various predictors based on
the notion of ARH process that is, \textit{linear interpolation ARH predictor}
(Pumo (1998)), \textit{Fourier interpolation ARF predictor} and \textit{ARW
predictor} based on the presentation (\ref{ARH1}), \ by two statistical
criteria: \textit{mean-squared error} (MSE) and \textit{relative mean-absolute
error} (RMAE) defined by :%
\[
MSE=\frac{1}{m}\sum_{j=1}^{m}(X_{n}(t_{j})-\hat{X}_{n}(t_{j}))^{2}%
,\hspace{0.11in}RMAE=\frac{1}{m}\sum_{j=1}^{m}\frac{|X_{n}(t_{j})-\hat{X}%
_{n}(t_{j})|}{|X_{n}(t_{j})|}%
\]
where $m$ is the number of discretized points.

The second example concerns real data, namely climatological time series
describing the El Ni\~{n}o-Southern Oscillation (ENSO). We compare our
predictor with predictor based on similar approaches found in the literature :
\textit{spline smoothing FAR predictor }(Besse and Cardot (1996)),
\textit{Local FAR predictor} (Besse et al. (2000))\ and \textit{wavelet based
predictor} (Antoniadis and Sapatinas (2003)).

\subsection{Some technical details about simulations}

Consider the Fourier basis on $L^{2}[0,\delta]$ and denote $e_{0}(t)=1/\delta$
and $e_{2j-1}(t)=\cos(2j\pi t/\delta)$, $e_{2j}(t)=\sin(2j\pi t/\delta
)$\ \ for $j\geq1$. Then a simple calculation shows that
\[
\mathbf{w}=\{e_{0},[1+4j^{2}\pi^{2}/\delta^{2}]^{-1/2}\cdot e_{2j-1}%
,[1+4j^{2}\pi^{2}/\delta^{2}]^{-1/2}\cdot e_{2j},j\geq1\}
\]
is an orthonormal basis for $W$. Let\ $f=%
%TCIMACRO{\tsum \nolimits_{i=0,\infty}}%
%BeginExpansion
{\textstyle\sum\nolimits_{i=0,\infty}}
%EndExpansion
c_{j}e_{j}$, where $c_{j}=\langle f,e_{j}\rangle_{L^{2}}$, be the Fourier
series of a continuously differentiable function $f$. Then $f^{\prime}=%
%TCIMACRO{\tsum \nolimits_{i=0,\infty}}%
%BeginExpansion
{\textstyle\sum\nolimits_{i=0,\infty}}
%EndExpansion
c_{j}e_{j}^{\prime}$. Furthermore the decomposition $%
%TCIMACRO{\tsum \nolimits_{j=0,\infty}}%
%BeginExpansion
{\textstyle\sum\nolimits_{j=0,\infty}}
%EndExpansion
\left\langle f,w_{j}\right\rangle _{W}w_{j}$ on $W$ of $f$ is given by $c_{0}+%
%TCIMACRO{\tsum \nolimits_{j=1,\infty}}%
%BeginExpansion
{\textstyle\sum\nolimits_{j=1,\infty}}
%EndExpansion
\left\langle f,w_{j}\right\rangle _{W}w_{j}$\ where $\left\langle
f,w_{2j-1}\right\rangle _{W}=[1+4j^{2}\pi^{2}/\delta^{2}]^{-1/2}%
\cdot\left\langle f,e_{2j-1}\right\rangle _{L^{2}}$ and $\left\langle
f,w_{2j}\right\rangle _{W}=[1+4j^{2}\pi^{2}/\delta^{2}]^{-1/2}\cdot
\left\langle f,e_{2j}\right\rangle _{L^{2}}$.

In order to calculate the covariance operators given in section
\ref{section.estimation.procedure} denote $\mathbf{w}_{N}$ (resp.
$\mathbf{e}_{N}$) $N$-vectors \ of the basis on $W$ (resp. $L^{2}$) that is
$\mathbf{w}_{N}=\,^{t}(w_{0},w_{1}\ldots w_{N-1})$\ (resp. $\mathbf{e}%
_{N}=\,^{t}(e_{0},e_{1},\ldots,e_{N-1})$) and suppose that $N$ is an odd and
positive number. Denote $\mathbf{X}$ and $\mathbf{X}^{\prime}$ the $N\times n$
matrices :%
\begin{align*}
\mathbf{X}  &  \mathbf{=}\ (\left\langle X_{i},w_{k-1}\right\rangle
_{W},\,k=1,\ldots,N;\,i=1,\ldots,n,\\
\mathbf{X}^{\prime}  &  \mathbf{=}\ (\left\langle X_{i}^{\prime}%
,e_{k-1}\right\rangle _{L^{2}},\,k=1,\ldots,N;\,i=1,\ldots,n.
\end{align*}
As noted above the coefficients $\left\langle X_{i},w_{k-1}\right\rangle _{W}$
and $\left\langle X_{i}^{\prime},e_{k-1}\right\rangle _{L^{2}}$ are obtained
directly from the Fourier decomposition of $X_{i}$, for $i=1,\ldots,n$.

It follows that covariance operators $\Gamma_{n},\Gamma_{n}^{\prime}%
,\Gamma_{n}^{\prime\ast},\Gamma_{n}^{^{\prime\prime}},\Delta_{n},\Delta
_{n}^{\prime}$\ can be approximated by :%
\begin{gather*}
C_{n}=(1/n)\cdot\mathbf{X\,(}^{t}\mathbf{X),\;}C_{n}^{\prime}=(1/n)\cdot
\mathbf{X}^{\prime}\mathbf{\,(}^{t}\mathbf{X),}\\
C_{n}^{\prime\ast}=(1/n)\cdot\mathbf{X\,(}^{t}\mathbf{X}^{\prime}%
\mathbf{),\;}C_{n}^{^{\prime\prime}}=(1/n)\cdot\mathbf{X}^{\prime}%
\mathbf{\,(}^{t}\mathbf{X}^{\prime}\mathbf{),}\\
D_{n}=(1/[n-1])\cdot\mathbf{X}_{-1}\mathbf{\,(}^{t}\mathbf{X}_{-n}%
\mathbf{),\;}D_{n}^{\prime}=(1/[n-1])\cdot\mathbf{X}_{-1}\mathbf{\,(}%
^{t}\mathbf{X}_{-n}^{\prime}\mathbf{)}%
\end{gather*}
where $\mathbf{X}_{-n}^{\prime}$\ (resp. $\mathbf{X}_{-1}$) is the matrix
$\mathbf{X}^{\prime}$ (resp. $\mathbf{X}$) without the column $n$ (resp. $1$).
So in order to obtain the estimators given in the second section it suffices
to substitute the covariance operators in (\ref{Snphi}-\ref{Tnpsi}) and
(\ref{estim})\textbf{ }by their approximations given above and choosing
suitable values for $\alpha_{n}$ and $\beta_{n}$.

\subsection{Wong process}

This process is defined for $u\in R$ by:%
\[
\xi_{u}=\sqrt{3}\,\exp\left(  -\sqrt{3}\,u\right)  \int_{0}^{\exp\left(
2u/\sqrt{3}\right)  }W_{s}\,ds.
\]
Wong process is a mean-square differentiable stationary Gaussian process which
is zero-mean and with variance $1$. Let $\delta>0$ and $X_{i}\in W$ given by
$X_{i+1}(t)=\xi_{i\cdot\delta+t}$ for $t\in]0,\delta]$. Let $\epsilon_{i+1}$
be a squared differentiable r.v. with values in $W^{2,1}[0,\delta]$ :
\begin{equation}
\epsilon_{i+1}(t)=\sqrt{3}\,\exp\left(  -\sqrt{3}\,(i+t)\right)  \int
_{\exp(2i/\sqrt{3})}^{\exp(2(i+t)/\sqrt{3})}(W_{s}-W_{\exp(2i/\sqrt{3})})\,ds.
\label{def.epsilon}%
\end{equation}
Then the process $(X_{i},i\in Z)$ can be written as:%

\begin{equation}
X_{i+1}=[\phi+\Psi(D)]X_{i}+\epsilon_{i+1} \label{ARHD.presentation}%
\end{equation}
where $c(t)=\frac{\sqrt{3}}{2}\cdot\exp(-\sqrt{3}t)\cdot\{\exp(2t/\sqrt
{3})-1\}$ and:%
\[
\lbrack\phi(f)](t)=[\exp(-\sqrt{3}t)+\sqrt{3}c(t)]f(1),\hspace{0.11in}%
[\Psi(D)(f)](t)=c(t)f^{\prime}(1).
\]
Furthermore $\epsilon_{i+1}$\ is independent of $X_{i},X_{i}^{\prime}$ and a
direct calculation shows that assumption $\mathbf{H1}$ is satisfied.

Using a method presented in Blanke and Pumo (2003) we simulated a Wong process
$\xi(t)$\ on $[0,192.65]$, that corresponds to $n=105$ intervals of length
$\delta=1.8348$, each known at $m=50$ equidistant points $t_{j},j=1,\ldots
,50$. A simulated process is presented in Figure \ref{wong_proc}. The
associated process $(X_{i})$ is a $W^{2,1}[0,1.8348]$\ valued process.

\begin{center}
\textbf{Please insert here Figure 1}
\end{center}

Fifty Wong processes were simulated and for each of them we calculated the MSE
and RMAE criteria. The mean values for the two criteria for the 50 simulations
and various predictors are presented in Table \ref{Wong_mean_table_errors}.
Figure \ref{wong_prediction} presents the different predictions for one of the simulations.

For the calculation of ARHD predictors we consider two values for the
parameter $\alpha_{n}$, that are 0.1 and 0.3. The corresponding values for
$\beta_{n}$ are 0.65 and 0.5. In the calculation of ARH, ARF and ARW predictor
we consider $k_{n}=1$, that is the projection subspace for observation is
equal to one (see Pumo (1998) for details). Simulations show that the ARW and
ARF predictor are very similar and when $m$ (this is the case for example when
$m=50$) is large they give similar results to linear interpolation ARH
predictor. But the three predictors are less better than the ARHD
predictor.\ Notice also that the choice of the optimal values for the
parameters $\alpha_{n}$ and $\beta_{n}$\ may be done by a cross-validation procedure.

\begin{center}
\textbf{Please insert here Table 1}

\textbf{Please insert here Figure 2}
\end{center}

\subsection{Example SST: Sea Surface Temperature}

The second example concerns a climatological time series describing the El
Ni\~{n}o-Southern Oscillation (see. for example Besse et al. (2000) or Smith
et al. (1996) for a description of the data\footnote{Data is freely avalilable
from http://www.cpc.ncep.noaa.gov/data/indices/index.html}). The series gives
the monthly mean El Ni\^{n}o sea surface temperature index from January 1950
to December 1969, that is $m=12$, and is presented in figure
\ref{sst_proc_until_96}. We compare the ARHD predictor with various functional
prediction methods.

\begin{center}
\textbf{Please insert here Figure 3}
\end{center}

In the first numerical study we compare the prediction of the temperature
during 1986 knowing the data until 1985. We calculated the ARHD predictor with
$\alpha_{n}=0.4$ and $0.1$ and $\beta_{n}=0.8$ and $0.4$. The MSE and RMAE
criteria for various functional predictors are given in Table
\ref{sst_prediction_1986_table}. Results show that the best method are Wavelet
II (one of the wavelet approaches proposed in Antoniadis and Sapatinas) and
spline smoothing FAR. Nevertheless our predictor is better than the other
predictors or the classical SARIMA $(0,1,1)\times(1,0,1)_{12}$ model (see for
example Brockwell and Davis (1987)). Figure \ref{sst_proc_until_96} displays
the observed data during 1986 and its predictors by some of the predictors
discussed above. Notice that the ARF or ARW predictors are not satisfactory as
$m=12$.

In the second numerical study we make 10 one year ahead forecasts for the
period 1986-97. The statistical criteria for various functional methods are
presented in Table \ref{sst_pred_87_96}. The reader may notice that the ARHD
method gives a similar prediction as the Local FAR method which is the best
functional prediction method appearing in Besse et al. (2000). Note finally
that as described in the introduction of this section the computational effort
to obtain an ARH predictor is comparable to that of the calculation of an ARH predictor.

\begin{center}
\textbf{Please insert here Table 2}

\textbf{Please insert here Figure 4}

\textbf{Please insert here Table 3}
\end{center}

\section{Proofs\label{section.proofs}}

\textbf{Proof of Proposition \ref{Ident}:}

The couple $\left(  \phi,\Psi\right)  $ will be identified whenever, for any
other couple $\left(  \widetilde{\phi},\widetilde{\Psi}\right)  $, if $\left(
\widetilde{\phi},\widetilde{\Psi}\right)  \Lambda=\left(  \phi,\Psi\right)
\Lambda,$ $\left(  \widetilde{\phi},\widetilde{\Psi}\right)  =\left(
\phi,\Psi\right)  $. This will be true if
\[
\left\{  \left(  U,V\right)  \in\mathcal{E}:\left(  U,V\right)  \Lambda
=0\right\}  =\left\{  0\right\}  .
\]

But $\Lambda$ may be decomposed as the product of three operators, namely :%
\[
\Lambda=\left(
\begin{array}
[c]{c}%
I\\
D
\end{array}
\right)  \Gamma\left(  I,D^{\ast}\right)  .
\]
As $\Gamma$ is one to one by assumption $\mathbf{H3},$ and since
\[
\left(  T,TD^{\ast}\right)  =0\mathrm{\ iff\ }T=0,
\]
it is readily seen that%
\[
\left(  U,V\right)  \Lambda=0\mathrm{\ iff\ }U+VD=0
\]
which finishes the proof of the Proposition.

We begin with five Lemmas needed to prove Theorem \ref{Conv}:

\begin{lemma}
\label{Bos}%
\begin{align*}
\left\Vert \Gamma_{n}-\Gamma\right\Vert _{\infty}  &  =O_{\mathbb{P}}\left(
\dfrac{1}{\sqrt{n}}\right)  ,\left\Vert \Gamma_{n}^{\prime}-\Gamma^{\prime
}\right\Vert _{\infty}=O_{\mathbb{P}}\left(  \dfrac{1}{\sqrt{n}}\right)  ,\\
\left\Vert \Gamma_{n}^{\prime\ast}-\Gamma^{\prime\ast}\right\Vert _{\infty}
&  =O_{\mathbb{P}}\left(  \dfrac{1}{\sqrt{n}}\right)  ,\left\Vert \Gamma
_{n}^{\prime\prime}-\Gamma^{\prime\prime}\right\Vert _{\infty}=O_{\mathbb{P}%
}\left(  \dfrac{1}{\sqrt{n}}\right)  .
\end{align*}

\end{lemma}

\textbf{Proof :} Since $X_{n}$ is an ARH(1) process (with autocorrelation
operator $A$), we can invoke for instance Theorem 4.1 p.98 in Bosq (2000) to
get the first rate of decay. All the other results above are due to the
boundedness (in our framework) of the differential operator $D.$ Indeed for
instance $\Gamma_{n}^{\prime\prime}=D\Gamma_{n}D^{\ast}$.

\begin{lemma}
\label{A1}%
\begin{align*}
\left\Vert \Gamma^{\prime\prime\dagger}\right\Vert _{\infty}  &  =\left\Vert
\Gamma_{n}^{\prime\prime\dagger}\right\Vert _{\infty}=\dfrac{1}{\alpha},\\
\left\Vert \Gamma^{\prime\prime\dagger}-\Gamma_{n}^{\prime\prime\dagger
}\right\Vert _{\infty}  &  =O_{\mathbb{P}}\left(  \dfrac{1}{\alpha^{2}\sqrt
{n}}\right)  .
\end{align*}

\end{lemma}

\textbf{Proof : }We prove the first bound%
\[
\left\Vert \Gamma^{\prime\prime\dagger}\right\Vert _{\infty}=\left\Vert
\left(  \Gamma^{\prime\prime}+\alpha I\right)  ^{-1}\right\Vert _{\infty}%
\]
and as $\Gamma^{\prime\prime}$ is a positive compact operator, the norm of
operator $\left(  \Gamma^{\prime\prime}+\alpha I\right)  ^{-1},$ which is
known as the resolvent operator of $\Gamma^{\prime\prime},$ is non random and
evaluated at $\alpha^{-1}.$ The same is true with $\Gamma_{n}^{\prime\prime}$
replacing $\Gamma^{\prime\prime}.$

Using $B^{-1}-A^{-1}=A^{-1}(A-B)B^{-1}$\ for two invertible matrices $A$ and
$B$, we get :%
\[
\Gamma^{\prime\prime\dagger}-\Gamma_{n}^{\prime\prime\dagger}=\left(
\Gamma_{n}^{\prime\prime}+\alpha I\right)  ^{-1}\left(  \Gamma_{n}%
^{\prime\prime}-\Gamma^{\prime\prime}\right)  \left(  \Gamma^{\prime\prime
}+\alpha I\right)  ^{-1}%
\]
which entails%
\begin{align*}
\left\Vert \Gamma^{\prime\prime\dagger}-\Gamma_{n}^{\prime\prime\dagger
}\right\Vert _{\infty}  &  \leq\left\Vert \left(  \Gamma_{n}^{\prime\prime
}+\alpha I\right)  ^{-1}\right\Vert _{\infty}\left\Vert \Gamma_{n}%
^{\prime\prime}-\Gamma^{\prime\prime}\right\Vert _{\infty}\left\Vert \left(
\Gamma^{\prime\prime}+\alpha I\right)  ^{-1}\right\Vert _{\infty}\\
&  =\dfrac{1}{\alpha^{2}}\left\Vert \Gamma_{n}^{\prime\prime}-\Gamma
^{\prime\prime}\right\Vert _{\infty}=O_{\mathbb{P}}\left(  \dfrac{1}%
{\alpha^{2}\sqrt{n}}\right)
\end{align*}
by Lemma \ref{Bos}.

\begin{lemma}
\label{A2}Let $S_{n,\phi}$ and $S_{\phi}$\ defined respectively by
(\ref{Snphi})\ and (\ref{Sphi}). Then%
\[
\left\Vert S_{n,\phi}-S_{\phi}\right\Vert _{\infty}=O_{\mathbb{P}}\left(
\dfrac{1}{\alpha^{2}\sqrt{n}}\right)  .
\]

\end{lemma}

\textbf{Proof }$:$ From equations (\ref{Sphi}) and (\ref{Spsi}) we obtain
$S_{n,\phi}-S_{\phi}=\Gamma_{n}-\Gamma+\Gamma^{\prime\ast}\left(
\Gamma^{\prime\prime\dagger}\right)  \Gamma^{\prime}-\Gamma_{n}^{\prime\ast
}\left(  \Gamma_{n}^{\prime\prime\dagger}\right)  \Gamma_{n}^{\prime}$ and%

\[
\left\Vert S_{n,\phi}-S_{\phi}\right\Vert _{\infty}\leq\left\Vert \Gamma
_{n}-\Gamma\right\Vert _{\infty}+\left\Vert \Gamma^{\prime\ast}\left(
\Gamma^{\prime\prime\dagger}\right)  \Gamma^{\prime}-\Gamma_{n}^{\prime\ast
}\left(  \Gamma_{n}^{\prime\prime\dagger}\right)  \Gamma_{n}^{\prime
}\right\Vert _{\infty}.
\]
We look for a bound for%
\begin{align*}
\lefteqn{\left\Vert \Gamma^{\prime\ast}\left(  \Gamma^{\prime\prime\dagger
}\right)  \Gamma^{\prime}-\Gamma_{n}^{\prime\ast}\left(  \Gamma_{n}%
^{\prime\prime\dagger}\right)  \Gamma_{n}^{\prime}\right\Vert _{\infty}%
\leq\left\Vert \Gamma^{\prime\ast}\left(  \Gamma^{\prime\prime\dagger}\right)
\Gamma^{\prime}-\Gamma_{n}^{\prime\ast}\left(  \Gamma^{\prime\prime\dagger
}\right)  \Gamma^{\prime}\right\Vert _{\infty}}\\
&  \hspace{0.59in}+\left\Vert \Gamma_{n}^{\prime\ast}\left(  \Gamma
^{\prime\prime\dagger}\right)  \Gamma^{\prime}-\Gamma_{n}^{\prime\ast}\left(
\Gamma^{\prime\prime\dagger}\right)  \Gamma_{n}^{\prime}\right\Vert _{\infty
}+\left\Vert \Gamma_{n}^{\prime\ast}\left(  \Gamma^{\prime\prime\dagger
}\right)  \Gamma_{n}^{\prime}-\Gamma_{n}^{\prime\ast}\left(  \Gamma
_{n}^{\prime\prime\dagger}\right)  \Gamma_{n}^{\prime}\right\Vert _{\infty}%
\end{align*}
Obviously the two first terms above may be bounded in probability by%
\[
\left\Vert \Gamma^{\prime\prime\dagger}\right\Vert _{\infty}\left\Vert
\Gamma_{n}^{\prime}\right\Vert _{\infty}\left\Vert \Gamma^{\prime}-\Gamma
_{n}^{\prime}\right\Vert _{\infty}=O_{\mathbb{P}}\left(  \dfrac{1}{\alpha
\sqrt{n}}\right)
\]
since $\left\Vert \Gamma_{n}^{\prime}\right\Vert _{\infty}=\left\Vert
\Gamma_{n}^{\prime\ast}\right\Vert _{\infty}.$ The remaining term may be
bounded by
\[
\left\Vert \Gamma_{n}^{\prime}\right\Vert ^{2}\left\Vert \Gamma^{\prime
\prime\dagger}-\Gamma_{n}^{\prime\prime\dagger}\right\Vert _{\infty
}=O_{\mathbb{P}}\left(  \left\Vert \Gamma^{\prime\prime\dagger}-\Gamma
_{n}^{\prime\prime\dagger}\right\Vert _{\infty}\right)
\]
and Lemma \ref{A1} finishes the proof.

\begin{center}

\end{center}

\begin{lemma}
\label{A2b}The operator $S_{\phi}$ is positive hence%
\[
\left\Vert \left(  S_{\phi}+\beta I\right)  ^{-1}\right\Vert _{\infty}%
\leq\dfrac{1}{\beta}.
\]

\end{lemma}

\textbf{Proof : }Before starting the proof it is worth reminding the following
fact to the reader. It $T$ is a compact operator from a Hilbert space $H_{1}$
to a Hilbert space $H_{2}$ it admits a Schmidt representation%
\[
T=\sum_{i=1}^{+\infty}s_{i}\left(  u_{i}\otimes v_{i}\right)
\]
where the $s_{i}$'s are the positive eigenvalues of $T^{\ast}T$ (i.e. of
$TT^{\ast}$) and where $u_{i}$ (resp. $v_{i}$) denote a complete orthonormal
system of $H_{1}$ (resp. $H_{2}$). We refer for instance to Theorem 1.1 page
96 in Gohberg, Goldberg, Kaashoek (1991). Now we turn to
\[
S_{\phi}=\Gamma-\Gamma^{\prime\ast}\Gamma^{\prime\prime\dagger}\Gamma^{\prime}%
\]
We set $U=D\Gamma^{1/2}$ (remind that $\Gamma^{\prime}=D\Gamma$). The operator
$U$ is a compact from $W$ to $L$ since $D$ is bounded and $\Gamma^{1/2}$ is
compact like $\Gamma.$ Then $\Gamma^{\prime\prime}=UU^{\ast}$ and
$\Gamma^{\prime\prime\dagger}=\left(  UU^{\ast}+\alpha I\right)  ^{-1}$ and we
rewrite%
\begin{equation}
S_{\phi}=\Gamma^{1/2}\left(  I-U^{\ast}\left(  UU^{\ast}+\alpha I\right)
^{-1}U\right)  \Gamma^{1/2} \label{posi}%
\end{equation}
Now let us write the Schmidt decomposition of $U$%
\[
U=\sum_{i=1}^{+\infty}s_{i}\left(  u_{i}\otimes_{W}v_{i}\right)  \quad
u_{i}\in W,v_{i}\in L
\]
Easy computations lead to
\[
\left(  I-U^{\ast}\left(  UU^{\ast}+\alpha I\right)  ^{-1}U\right)
=\sum_{i=1}^{+\infty}\dfrac{\alpha}{\alpha+s_{i}^{2}}u_{i}\otimes_{W}u_{i}%
\]
From (\ref{posi}) we deduce that for all $x$ in $W$ $\left\langle S_{\phi
}x,x\right\rangle \geq0$ hence the announced result.

\begin{lemma}
\label{A3}
\[
\left\Vert \left(  S_{n,\phi}+\beta I\right)  ^{-1}-\left(  S_{\phi}+\beta
I\right)  ^{-1}\right\Vert _{\infty}=O_{\mathbb{P}}\left(  \dfrac{1}%
{\alpha^{2}\beta^{2}\sqrt{n}}\right)  .
\]

\end{lemma}

\textbf{Proof :}%
\begin{equation}
\left(  S_{n,\phi}+\beta I\right)  ^{-1}-\left(  S_{\phi}+\beta I\right)
^{-1}=\left(  S_{\phi}+\beta I\right)  ^{-1}\left(  S_{\phi}-S_{n,\phi
}\right)  \left(  S_{n,\phi}+\beta I\right)  ^{-1} \label{bec}%
\end{equation}
hence%
\begin{equation}
\left(  I-\left(  S_{\phi}+\beta I\right)  ^{-1}\left(  S_{\phi}-S_{n,\phi
}\right)  \right)  \left(  S_{n,\phi}+\beta I\right)  ^{-1}=\left(  S_{\phi
}+\beta I\right)  ^{-1}. \label{diese}%
\end{equation}
Since
\[
\left\Vert \left(  S_{\phi}+\beta I\right)  ^{-1}\left(  S_{\phi}-S_{n,\phi
}\right)  \right\Vert _{\infty}\overset{_{\mathbb{P}}}{\rightarrow}0
\]
the probability that $I-\left(  S_{\phi}+\beta I\right)  ^{-1}\left(  S_{\phi
}-S_{n,\phi}\right)  $ is an invertible operator tends to 1. It suffices
indeed that%
\[
\left\Vert \left(  S_{\phi}+\beta I\right)  ^{-1}\left(  S_{\phi}-S_{n,\phi
}\right)  \right\Vert _{\infty}<1
\]
to write from (\ref{diese}) :%
\begin{equation}
\left(  S_{n,\phi}+\beta I\right)  ^{-1}=\left(  I-\left(  S_{\phi}+\beta
I\right)  ^{-1}\left(  S_{\phi}-S_{n,\phi}\right)  \right)  ^{-1}\left(
S_{\phi}+\beta I\right)  ^{-1}. \label{bem}%
\end{equation}
We set $H_{n}=\left(  S_{\phi}+\beta I\right)  ^{-1}\left(  S_{\phi}%
-S_{n,\phi}\right)  $ then%
\[
\left(  I-H_{n}\right)  ^{-1}=I+\sum_{p=1}^{+\infty}H_{n}^{p}.
\]
At last from (\ref{bem}) and (\ref{bec}) :%
\[
\left(  S_{n,\phi}+\beta I\right)  ^{-1}-\left(  S_{\phi}+\beta I\right)
^{-1}=\sum_{p=1}^{+\infty}H_{n}^{p}\left(  S_{\phi}+\beta I\right)  ^{-1}%
\]
and
\begin{align*}
\left\Vert \left(  S_{n,\phi}+\beta I\right)  ^{-1}-\left(  S_{\phi}+\beta
I\right)  ^{-1}\right\Vert _{\infty}  &  \leq\left\Vert \left(  S_{\phi}+\beta
I\right)  ^{-1}\right\Vert _{\infty}\sum_{p=1}^{+\infty}\left\Vert
H_{n}\right\Vert _{\infty}\\
&  =O_{\mathbb{P}}\left(  \dfrac{1}{\alpha^{2}\beta^{2}\sqrt{n}}\right)
\end{align*}
by Lemmas \ref{A2} and \ref{A2b}.

\textbf{Proof of Theorem \ref{Conv} : }We prove the Proposition for $\phi_{n}$
since the same technique would lead to an analogous result for $\Psi$.
Developing the expression of $\Delta_{n}$ and $\Delta_{n}^{\prime}$ yields%
\begin{align*}
\Delta_{n}  &  =\phi\Gamma_{n}+\Psi\Gamma_{n}^{\prime}+U_{n},\\
\Delta_{n}^{\prime}  &  =\phi\Gamma_{n}^{\prime\ast}+\Psi\Gamma_{n}%
^{\prime\prime}+U_{n}^{\prime},
\end{align*}
with
\[
U_{n}=\dfrac{1}{n}\sum_{k=1}^{n}X_{k}\otimes\varepsilon_{k+1},\;U_{n}^{\prime
}=\dfrac{1}{n}\sum_{k=1}^{n}X_{k}^{\prime}\otimes\varepsilon_{k+1}.
\]
Hence
\begin{align*}
T_{n,\phi}  &  =\Delta_{n}-\Delta_{n}^{\prime}\left(  \Gamma_{n}^{\prime
\prime\dagger}\right)  \Gamma_{n}^{\prime}\\
&  =\phi\Gamma_{n}+\Psi\Gamma_{n}^{\prime}-\phi\Gamma_{n}^{\prime\ast}\left(
\Gamma_{n}^{\prime\prime\dagger}\right)  \Gamma_{n}^{\prime}-\Psi\Gamma
_{n}^{\prime\prime}\left(  \Gamma_{n}^{\prime\prime\dagger}\right)  \Gamma
_{n}^{\prime}+U_{n}-U_{n}^{\prime}\left(  \Gamma_{n}^{\prime\prime\dagger
}\right)  \Gamma_{n}^{\prime}\\
&  =\phi S_{n,\phi}+\Psi\left[  \Gamma_{n}^{\prime}-\Gamma_{n}^{\prime\prime
}\left(  \Gamma_{n}^{\prime\prime\dagger}\right)  \Gamma_{n}^{\prime}\right]
+\left[  U_{n}-U_{n}^{\prime}\left(  \Gamma_{n}^{\prime\prime\dagger}\right)
\Gamma_{n}^{\prime}\right]
\end{align*}
At last%
\begin{align}
\phi_{n}-\phi &  =\beta\phi\left(  S_{n,\phi}+\beta I\right)  ^{-1}%
+\Psi\left[  \Gamma_{n}^{\prime}-\Gamma_{n}^{\prime\prime}\left(  \Gamma
_{n}^{\prime\prime\dagger}\right)  \Gamma_{n}^{\prime}\right]  \left(
S_{n,\phi}+\beta I\right)  ^{-1}\label{dec}\\
&  \hspace{0.37in}+\left[  U_{n}-U_{n}^{\prime}\left(  \Gamma_{n}%
^{\prime\prime\dagger}\right)  \Gamma_{n}^{\prime}\right]  \left(  S_{n,\phi
}+\beta I\right)  ^{-1}.\nonumber
\end{align}
The proof will be achieved if we prove that the three terms in the display
above tend to zero in probability. The three next Propositions namely
Propositions \ref{T3}, \ref{T1} and \ref{T2} are devoted to this goal. We
begin with the last one involving $U_{n}$ and $U_{n}^{\prime}.$ We need two
auxiliary Lemmas.

\begin{lemma}
\label{Un}%
\begin{align}
\left\Vert U_{n}\right\Vert  &  =O_{\mathbb{P}}\left(  \dfrac{1}{\sqrt{n}%
}\right)  ,\label{G1}\\
\left\Vert U_{n}^{\prime}\right\Vert  &  =O_{\mathbb{P}}\left(  \dfrac
{1}{\sqrt{n}}\right)  ,\label{G2}\\
\left\Vert U_{n}^{\prime}\left(  \Gamma_{n}^{\prime\prime\dagger}\right)
\Gamma_{n}^{\prime}\right\Vert _{\infty}  &  =O_{\mathbb{P}}\left(  \dfrac
{1}{\alpha\sqrt{n}}\right)  . \label{G3}%
\end{align}

\end{lemma}

\textbf{Proof : }The proof of (\ref{G1}) and (\ref{G2}) is obvious since
$U_{n}$ and $U_{n}^{\prime}$ are sums of uncorrelated random operators (here
uncorrelated means that the cross covariance operator between two distinct
random elements is the null operator). Then%
\[
\left\Vert U_{n}^{\prime}\left(  \Gamma_{n}^{\prime\prime\dagger}\right)
\Gamma_{n}^{\prime}\right\Vert _{\infty}\leq\left\Vert U_{n}^{\prime
}\right\Vert _{\infty}\left\Vert \Gamma_{n}^{\prime\prime\dagger}\right\Vert
_{\infty}\left\Vert \Gamma_{n}^{\prime}\right\Vert _{\infty}%
\]
where the last term on the right side is bounded in probability, the first is
an $O_{\mathbb{P}}\left(  n^{-1/2}\right)  $ and the norm of the second is
almost surely bounded by $\alpha^{-1}$ and (\ref{G3}) is proved.

\begin{proposition}
\label{T3}%
\[
\left[  U_{n}-U_{n}^{\prime}\left(  \Gamma_{n}^{\prime\prime\dagger}\right)
\Gamma_{n}^{\prime}\right]  \left(  S_{n,\phi}+\beta I\right)  ^{-1}%
=O_{\mathbb{P}}\left(  \dfrac{1}{\alpha\beta\sqrt{n}}\right)
\]

\end{proposition}

\textbf{Proof : }The proof of the Proposition is a consequence of Lemmas
\ref{A2b} and \ref{Un}.

We turn to the first term in (\ref{dec}).

\begin{proposition}
\label{T1}If $\beta\rightarrow0$ and $\dfrac{1}{\alpha^{2}\beta^{2}\sqrt{n}%
}\rightarrow0,$%
\[
\beta\phi\left(  S_{n,\phi}+\beta I\right)  ^{-1}\overset{\mathbb{P}%
}{\rightarrow}0
\]

\end{proposition}

\textbf{Proof : }We invoke Lemma \ref{A3} to claim that it suffices to drop
the index $n$ in the Proposition and to prove that :%
\[
\beta\phi\left(  S_{\phi}+\beta I\right)  ^{-1}\overset{\mathbb{P}%
}{\rightarrow}0.
\]
In fact Lemma \ref{A3} links the asymptotic behavior of $\left(  S_{n,\phi
}+\beta I\right)  ^{-1}$ and $\left(  S_{\phi}+\beta I\right)  ^{-1}$.
Remember that $\phi$ is a compact operator from $W$ to $W.$ This fact is
crucial. It implies that we just have to prove that
\[
\beta\left(  S_{\phi}+\beta I\right)  ^{-1}x\rightarrow0
\]
By Lemma \ref{A2b} this fact is straightforward. Indeed it was then proved
that $S_{\phi}$ is a selfadjoint positive operator hence admits the spectral
decomposition%
\[
S_{\phi}=\sum_{i=1}^{+\infty}\mu_{i}\left(  t_{i}\otimes_{W}t_{i}\right)
\]
where the $\mu_{i}$'s are the positive eigenvalues of $S_{\phi}$ arranged in a
decreasing order and the $t_{i}$'s are the associated eigenvectors. Then if
$x=\sum x_{i}t_{i}$ where $x_{i}=\left\langle x,t_{i}\right\rangle _{W}$ we
easily get%
\[
\left\Vert \beta\left(  S_{\phi}+\beta I\right)  ^{-1}x\right\Vert _{W}%
^{2}=\sum_{i=1}^{+\infty}\dfrac{\beta^{2}\left\langle x,t_{i}\right\rangle
_{W}^{2}}{\left(  \mu_{i}+\beta\right)  ^{2}}%
\]
For fixed $i$, $\dfrac{\beta^{2}}{\left(  \mu_{i}+\beta\right)  ^{2}%
}\rightarrow0$ as $\beta\rightarrow0$ , $\sup_{i}\left\{  \dfrac{\beta^{2}%
}{\left(  \mu_{i}+\beta\right)  ^{2}}\right\}  \leq1$ and as $\sum
_{i=1}^{+\infty}\left\langle x,t_{i}\right\rangle _{W}^{2}<+\infty,$ applying
Lebesgue's dominated convergence Theorem yields $\beta\left(  S_{\phi}+\beta
I\right)  ^{-1}x\rightarrow0$ in $W.$\bigskip

Let us deal with the second term in (\ref{dec}).

\begin{proposition}
\label{T2}if $\beta$ and $\sqrt{\alpha}/\beta$ decays to zero,%
\[
\left\Vert \Psi\left[  \Gamma_{n}^{\prime}-\Gamma_{n}^{\prime\prime}\left(
\Gamma_{n}^{\prime\prime\dagger}\right)  \Gamma_{n}^{\prime}\right]  \left(
S_{n,\phi}+\beta I\right)  ^{-1}\right\Vert _{\infty}\overset{\mathbb{P}%
}{\rightarrow}0.
\]

\end{proposition}

\textbf{Proof : }Once more it suffices to prove that
\[
\dfrac{\left\Vert \left[  \Gamma_{n}^{\prime}-\Gamma_{n}^{\prime\prime}\left(
\Gamma_{n}^{\prime\prime\dagger}\right)  \Gamma_{n}^{\prime}\right]
\right\Vert _{\infty}}{\beta}\overset{\mathbb{P}}{\rightarrow}0.
\]
since%
\begin{equation}
\left\Vert \left[  \Gamma_{n}^{\prime}-\Gamma_{n}^{\prime\prime}\left(
\Gamma_{n}^{\prime\prime\dagger}\right)  \Gamma_{n}^{\prime}\right]  \left(
S_{n,\phi}+\beta I\right)  ^{-1}\right\Vert _{\infty}\leq\left(
\dfrac{\left\Vert \left[  \Gamma_{n}^{\prime}-\Gamma_{n}^{\prime\prime}\left(
\Gamma_{n}^{\prime\prime\dagger}\right)  \Gamma_{n}^{\prime}\right]
\right\Vert _{\infty}}{\beta}\right)  . \label{zeend}%
\end{equation}
We keep on replacing the random operators based on the sample by their limits.%
\begin{align*}
\lefteqn{\left\Vert \left[  \Gamma_{n}^{\prime}-\Gamma_{n}^{\prime\prime
}\left(  \Gamma_{n}^{\prime\prime\dagger}\right)  \Gamma_{n}^{\prime}\right]
\right\Vert _{\infty}=\left\Vert \alpha\left(  \Gamma_{n}^{\prime\prime
\dagger}\right)  \Gamma_{n}^{\prime}\right\Vert _{\infty}}  &  & \\
&  & \hspace{0.5in} \leq\left\Vert \alpha\left(  \Gamma_{n}^{\prime
\prime\dagger}\right)  \right\Vert _{\infty}\left\Vert \Gamma_{n}^{\prime
}-\Gamma^{\prime}\right\Vert _{\infty}+\left\Vert \alpha\left(  \Gamma
_{n}^{\prime\prime\dagger}-\Gamma^{\prime\prime\dagger}\right)  \Gamma
^{\prime}\right\Vert _{\infty}+\left\Vert \alpha\left(  \Gamma^{\prime
\prime\dagger}\right)  \Gamma^{\prime}\right\Vert _{\infty}%
\end{align*}

By Lemma \ref{Bos} the first term is an $O_{\mathbb{P}}\left(  \dfrac{1}%
{\sqrt{n}}\right)  ,$ the second is an $O_{\mathbb{P}}\left(  \dfrac{1}%
{\alpha\sqrt{n}}\right)  $.

The last term is totally deterministic and we are going to prove that it is an
$o\left(  \sqrt{\alpha}\right)  $. Once again we introduce the compact
operator $U=D\Gamma^{1/2}.$ We see that
\[
\alpha\Gamma^{\prime\prime\dagger}\Gamma^{\prime}=\alpha\left(  \Gamma
^{\prime\prime\dagger}\right)  ^{1/2}\left(  UU^{\ast}+\alpha I\right)
^{-1/2}U\Gamma^{1/2}%
\]
since $\Gamma^{\prime\prime}=UU^{\ast}.$ First we show that $\left(  UU^{\ast
}+\alpha I\right)  ^{-1/2}U$ is a class of operators uniformely bonded with
respect to $\alpha.$ In fact introducing the Schmidt representation for $U$
from Lemma \ref{A2b} we get
\[
\left(  UU^{\ast}+\alpha I\right)  ^{-1/2}U=\sum_{i=1}^{+\infty}\dfrac{s_{i}%
}{\sqrt{s_{i}^{2}+\alpha}}u_{i}\otimes u_{i}%
\]
and
\[
\sup_{\alpha\geq0}\left\Vert \left(  UU^{\ast}+\alpha I\right)  ^{-1/2}%
U\right\Vert _{\infty}\leq\sup_{\alpha,i}\left\{  \dfrac{s_{i}}{\sqrt
{s_{i}^{2}+\alpha}}\right\}  \leq1
\]
At last noting that $\left\Vert \alpha\left(  \Gamma^{\prime\prime\dagger
}\right)  ^{1/2}\right\Vert _{\infty}=o\left(  \sqrt{\alpha}\right)  $ and
taking into account (\ref{zeend}), the proof of Proposition \ref{T2} is
finished.

%

%TCIMACRO{\FRAME{ftbpFU}{6.6608in}{4.062in}{0pt}{\Qcb{Wong process simulated on
%$[0,192.65]$}}{\Qlb{wong_proc}}{wong_proc.wmf}%
%{\special{ language "Scientific Word";  type "GRAPHIC";  display "USEDEF";
%valid_file "F";  width 6.6608in;  height 4.062in;  depth 0pt;
%original-width 11.6931in;  original-height 8.2659in;  cropleft "0";
%filename '../Figures/wong_proc.WMF';file-properties "XNPEU";}}}%
%BeginExpansion
\begin{figure}
[ptb]
\begin{center}
\includegraphics[
natheight=8.265900in,
natwidth=11.693100in,
height=4.062in,
width=6.6608in
]%
{../Figures/wong_proc.wmf}%
\caption{Wong process simulated on $[0,192.65]$}%
\label{wong_proc}%
\end{center}
\end{figure}
%EndExpansion

\begin{table}[ptb]
\begin{center}%
\begin{tabular}
[c]{|c|c|c|c|c|c|}\hline
& \multicolumn{3}{|c}{ARH approach} & \multicolumn{2}{|c|}{ARHD approach}%
\\\hline
& $ARH$ & $ARF$ & $ARW$ & $\alpha=0.3,\beta=0.65$ & $\alpha=0.1,\beta
=0.5$\\\hline
$MSE$ & 0.624 & 0.622 & 0.623 & 0.327 & 0.323\\
$RMAE$ & 1.580 & 1.599 & 1.599 & 1.223 & 1.125\\\hline
\end{tabular}
\end{center}
\caption{Mean of MSE and RMAE error for the $50$ simulations.}%
\label{Wong_mean_table_errors}%
\end{table}

\bigskip%
%TCIMACRO{\FRAME{ftbpFU}{5.5357in}{3.9202in}{0pt}{\Qcb{Prediction of the
%$105$-th sample path by $ARH$, $ARF$, $ARW$ and $ARHD$ method.}}%
%{\Qlb{wong_prediction}}{wong_prediction.wmf}%
%width 5.5357in;  height 3.9202in;  depth 0pt;  original-width 11.6931in;
%original-height 8.2659in;  cropleft "0";  croptop "1";  cropright "1";
%filename 'Figures/wong_prediction.WMF';file-properties "XNPEU";}}}%
%BeginExpansion
\begin{figure}
[ptb]
\begin{center}
\includegraphics[
natheight=8.265900in,
natwidth=11.693100in,
height=3.9202in,
width=5.5357in
]%
{Figures/wong_prediction.wmf}%
\caption{Prediction of the $105$-th sample path by $ARH$, $ARF$, $ARW$ and
$ARHD$ method.}%
\label{wong_prediction}%
\end{center}
\end{figure}
%EndExpansion

\bigskip%
%TCIMACRO{\FRAME{ftbpFU}{5.8582in}{3.4835in}{0pt}{\Qcb{The monthly mean
%Ni\^{n}o-3 sea surface temperature index from 1950 until 1996.}}%
%{\Qlb{sst_proc_until_96}}{sst_proc_until_96.wmf}%
%{\special{ language "Scientific Word";  type "GRAPHIC";  display "USEDEF";
%valid_file "F";  width 5.8582in;  height 3.4835in;  depth 0pt;
%original-width 11.6931in;  original-height 8.2659in;  cropleft "0";
%filename '../Figures/sst_proc_until_96.WMF';file-properties "XNPEU";}}}%
%BeginExpansion
\begin{figure}
[ptb]
\begin{center}
\includegraphics[
natheight=8.265900in,
natwidth=11.693100in,
height=3.4835in,
width=5.8582in
]%
{../Figures/sst_proc_until_96.wmf}%
\caption{The monthly mean Ni\^{n}o-3 sea surface temperature index from 1950
until 1996.}%
\label{sst_proc_until_96}%
\end{center}
\end{figure}
%EndExpansion

\begin{table}[ptb]
\begin{center}%
\begin{tabular}
[c]{|l|c|c|}\hline
\multicolumn{3}{|c|}{Error of predictions 86}\\\hline
Predictor & MSE & RMAE\\\hline
Wavelet II & 0.063 & 0.89\%\\
FAR & 0.065 & 0,89\%\\
ARHD $\alpha=0.1,\beta=0.4$ & 0.167 & 1.25\%\\
Wavelet III & 0.191 & 1.20\%\\
ARHD $\alpha=0.4,\beta=0.8$ & 0.219 & 1.33\%\\
ARH(1) $k_{n}=1$ & 0.278 & 1.60\%\\
SARIMA & 1.457 & 3.72\%\\\hline
\end{tabular}
\end{center}
\caption{MSE and RMAE for the prediction of El Ni\^{n}o surface temperatures
during 1986 for various methods. }%
\label{sst_prediction_1986_table}%
\end{table}\bigskip%
%TCIMACRO{\FRAME{ftbpFU}{5.7726in}{4.0888in}{0pt}{\Qcb{Evolution of Ni\^{n}o-3
%surface temperature during 1986 and its various predictions.}}%
%{\Qlb{sst_prediction_86_fig}}{sst_prediction_86.wmf}%
%width 5.7726in;  height 4.0888in;  depth 0pt;  original-width 11.6931in;
%original-height 8.2659in;  cropleft "0";  croptop "1";  cropright "1";
%filename 'Figures/sst_prediction_86.WMF';file-properties "XNPEU";}}}%
%BeginExpansion
\begin{figure}
[ptbptb]
\begin{center}
\includegraphics[
natheight=8.265900in,
natwidth=11.693100in,
height=4.0888in,
width=5.7726in
]%
{Figures/sst_prediction_86.wmf}%
\caption{Evolution of Ni\^{n}o-3 surface temperature during 1986 and its
various predictions.}%
\label{sst_prediction_86_fig}%
\end{center}
\end{figure}
%EndExpansion

\begin{table}[ptb]
\begin{center}%
\begin{tabular}
[c]{|l|c|c|}\hline
\multicolumn{3}{|c|}{Mean Error of predictions 87-96}\\\hline
Predictor & MSE & RMAE\\\hline
ARHD $\alpha=0.4,\beta=0.8$ & 0.53 & 2.1\%\\
Local FAR & 0.53 & 2.2\%\\
ARHD $\alpha=0.1,\beta=0.4$ & 0.53 & 2.2\%\\
FAR & 0.55 & 2.3\%\\
ARH(1) $k_{n}=1$ & 0.68 & 2.4\%\\
SARIMA & 1.45 & 3.7\%\\\hline
\end{tabular}
\end{center}
\caption{Mean value of MSE and RMAE errors for prediction of SST from 1987 to
1996.}%
\label{sst_pred_87_96}%
\end{table}
\end{document}